\documentclass[letterpaper, 10 pt, conference]{ieeeconf}  

\IEEEoverridecommandlockouts                 
\usepackage{multibib}
\usepackage{graphicx}
\usepackage{pst-all}
\usepackage{amsmath}
\usepackage{hyperref}

\title{\LARGE \bf
3D Printing \& the Artists Known As the Fauves \\
\textit{July 2016} 
}

\author{Prof. Edward Aboufadel\\
Grand Valley State University,\\
Allendale, MI, USA\\
aboufade@gvsu.edu
}

\begin{document}

\maketitle
\thispagestyle{empty}
\pagestyle{empty}

\begin{abstract}

In this short, chatty paper, I describe how my attempt to use mathematics to create a 3D print of a school portrait led me a group of early 20\textsuperscript{th} century French artists known as the Fauves.

\end{abstract}

\section{3D Printing:  Something From Nothing}

At the beginning of the decade, I started to hear about low-cost 3D printers, and I decided that I wanted one.  As the end of the 2012-2013 fiscal year approached, I asked my Dean to purchase a Makerbot Replicator 2 for my department, and he agreed.  Soon, we had a new 3D printer set up in my office, and I had no idea what to do with it.  With two REU students that summer, we started asking the question, "What can a mathematician do with a 3D printer?"

The first thing to do was to ``print'' something.  The Makerbot came with several stereolithography (STL) files, which are large files that describe the object as a set of oriented triangles.  STL files are processed by software known as a ``slicer'', for reasons that will be obvious in a moment, and the results of the slicer direct the printer's actions.  Following these directions, a 3D printer creates an object by melting filament (a spaghetti-like spool of plastic, usually a type known as Polylactic acid, or PLA) at temperatures over 200 degrees Celsius, and then extruding the hot plastic through a nozzle onto a build plate, layer-by-layer, or slice-by-slice.

\begin{figure}[thpb]
      \centering
      \includegraphics[scale=1.5]{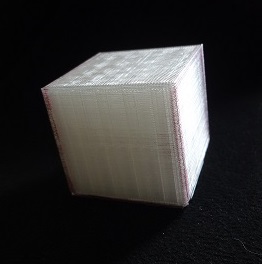}
      \caption{3D print of a cube.}
      \label{cube}
\end{figure}

Mathematically, working in three dimensions, the printer starts by laying down plastic in the plane $z=0$, at various points $(x, y, 0)$, following paths created by the slicer \cite{Austin14}.  Then the build plate moves down vertically a very small amount $\Delta z$ so that more plastic can be extruded in the plane $z = \Delta z$.  (The nozzle moves in the $x$ and $y$ directions, but not the $z$ direction.)  The material for this second slice usually is placed on top of the material in the $z=0$ plane.  The printer continues to $z = 2 \Delta z$, $z = 3 \Delta z$, etc., and there are a variety of rules to take into account curved surfaces, overhangs, and other aspects of a three-dimensional object, or else all we could ever print would be rectangular prisms and pyramids.  One of our first prints can be found in Figure \ref{cube}.

\begin{figure}[thpb]
        \centering
        \includegraphics[scale=1.25]{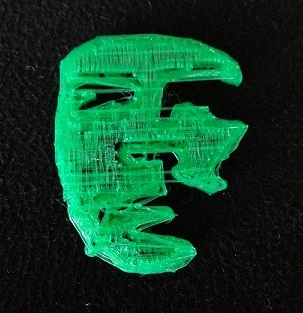}
        \vspace{1cm}
        \includegraphics[scale=0.2]{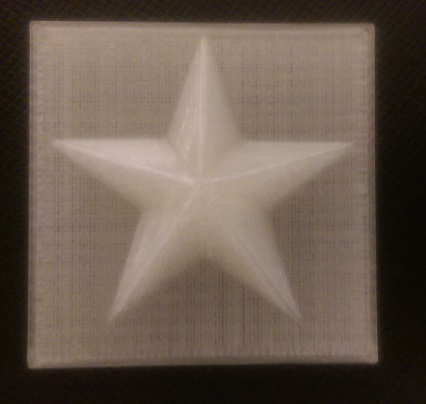}
        \vspace{1cm}
        \includegraphics[scale=1.5]{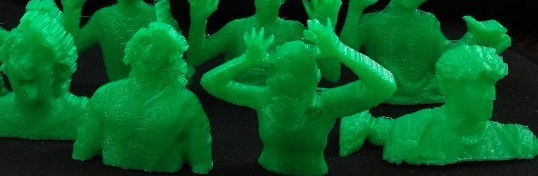}
        \vspace{1cm}
        \includegraphics[scale=0.08]{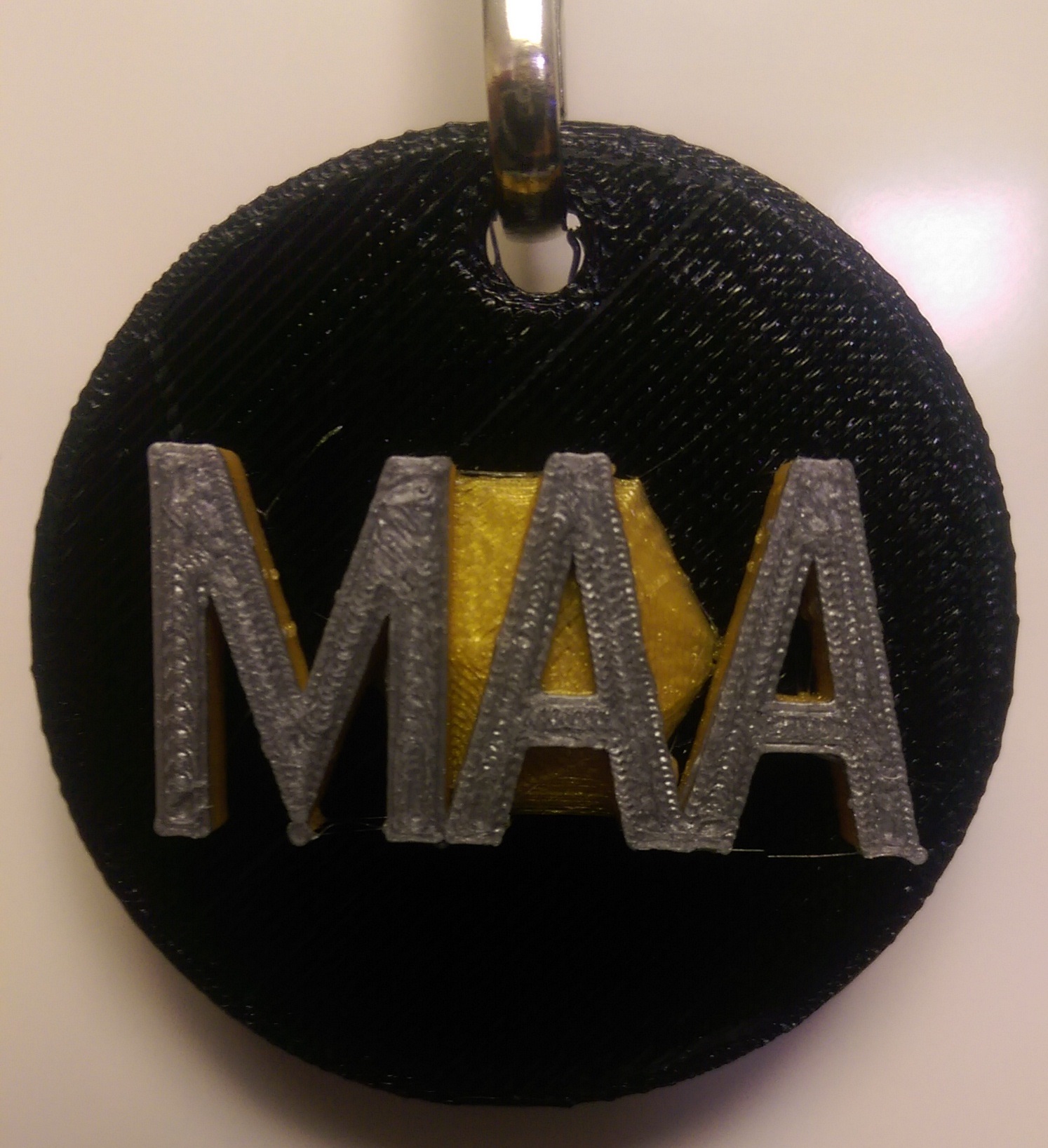}
		\caption{3D printing examples:  Che; Star; Data capture from Kinect; Pendant celebrating the Mathematical Association of America's 100\textsuperscript{th} anniversary.}
        \label{samples}
\end{figure}

\begin{figure}[thpb]
      \centering
       \includegraphics[scale=0.3]{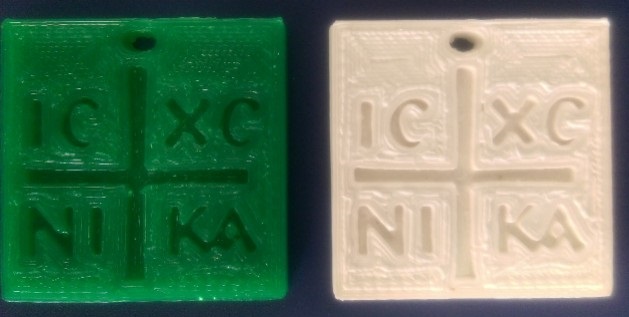}
      \caption{Victor's cross Christmas ornament.}
      \label{victors}
   \end{figure}

There are a variety of ways to create STL files.  For very simple designs, STL files can be written in a text editor, where the vertices and outer normal vector of each triangle are listed, such as:

\begin{verbatim}
    facet normal 0.0 1.0 0.0
        outer loop
          vertex 0.0 40.0 0.0
          vertex 40.0 40.0 0.0
          vertex 0.0 40.0 40.0
        endloop
    endfacet
\end{verbatim}

However, nearly all STL files are created with a software program of some sort, and saved in relatively efficient binary files.  These files can be created using 3D analogues of Microsoft Paint available for free online such as Tinkercad (\url{https://www.tinkercad.com/}), open source programs like OpenScad (\url{www.openscad.org}), and proprietary programs like \textit{Mathematica}.  Influenced by Segerman's paper \cite{Segerman}, as well as my affinity for mathematical analysis (as opposed to geometry), my students and I have primarily used \textit{Mathematica} to design objects for 3D printing.  Some of our creations from the past two years can be seen in Figures \ref{samples} and \ref{victors}.

\section{Captain, I'm An Analyst, Not A Geometer}

There are other mathematicians who have designed and printed fantastic objects with 3D printers in the past four years.  Laura Taalman spent a year doing daily mathematically-inspired 3D-printing projects that she documented on her blog (\url{http://makerhome.blogspot.com/}).  George Hart has created 3D printed fractals (\url{https://www.simonsfoundation.org/multimedia/3-d-printing-of-mathematical-models/}).  Links to many other projects can be found on my web site (\url{https://sites.google.com/site/aboufadelreu/Profile/3d-printing}).  Most of these projects are based on mathematical topics like knot theory, polyhedra, and minimal surfaces that have geometry as a foundation.

However, I have always been more of an analyst, thinking in terms of functions defined on a domain, so my 3D printing projects have had more of a flavor of computing a function $f(x,y)$ that represents a height above points on a compact set in the $xy$-plane.  The MAA pendant in Figure 2 is a good example of this approach, with the domain being a circle of the form $(x-h)^2 + (y - k)^2 = r^2$, and the function having a value of 0 where the hole for the clip is, 2 (millimeters, which is the standard unit in 3D printing design) for the parts of the outer parts of the pendant, and then something significantly more elaborate for the combination of the initials and the icosohedron.

\begin{figure}[thpb]
      \centering
      \includegraphics[scale=0.16]{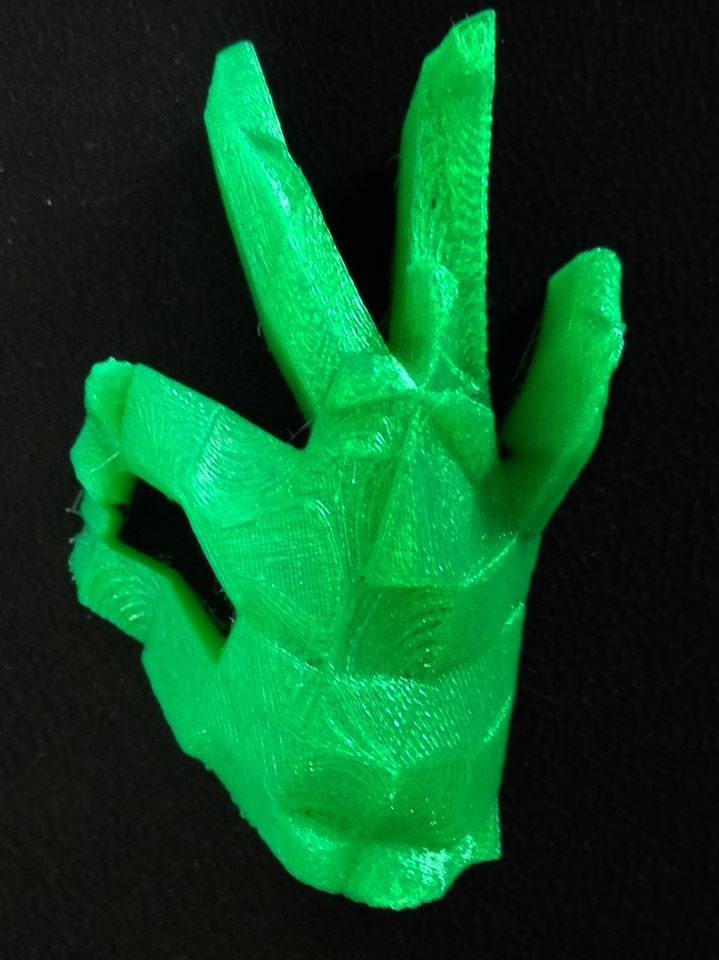}
      \caption{Based on two photographs of a hand.}
      \label{handprint}
   \end{figure}

The first major project with my REU students was to use two photographs of the same object (specifically, a hand) to create a 3D print of the object, as seen in Figure \ref{handprint}.  The idea was that the stereo effect of two photographs, like the way our pair of eyes works, allows us to deduce the distance of various points on the object to the cameras, and from there a depth function can be defined \cite{hand}.  After we completed this project, I started to wonder to what extent this could be done with just one photograph.  This is not a new question (for example, there is the new ``123D Catch'' app from Autodesk, see \url{http://www.123dapp.com/catch}).  However, I wanted to figure out how to do this with my methods:  create and generate a 3D rendering of a function $f(x,y)$ in \textit{Mathematica} that represents height or depth, export the rendering as an STL file, then slice and print.  Of course, the hard part is the creation of the function.  And my naive use of grayscale, which I will describe in the next section, wasn't going to get me very far.

\section{Fifty Shades of Grayscale}

One of the tricks to create $f(x,y)$ from an image like a logo \cite{logo} is to base the values of the function on the shades in a grayscale image.  Grayscale describes intensity (or, as we will discover below, \textit{luminance}), with brighter shades having larger grayscale values.  There are two equivalent versions of grayscale:  a scale of 0 (black) to 1 (white), and a scale of 0 to 255 (which comes from using bit strings).

\begin{figure}[thpb]
      \centering
       \includegraphics[scale=1]{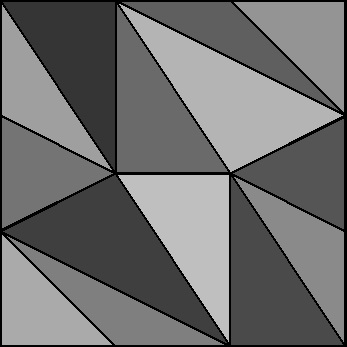}
       
       \vspace{1cm}
       
      \includegraphics[scale=0.08]{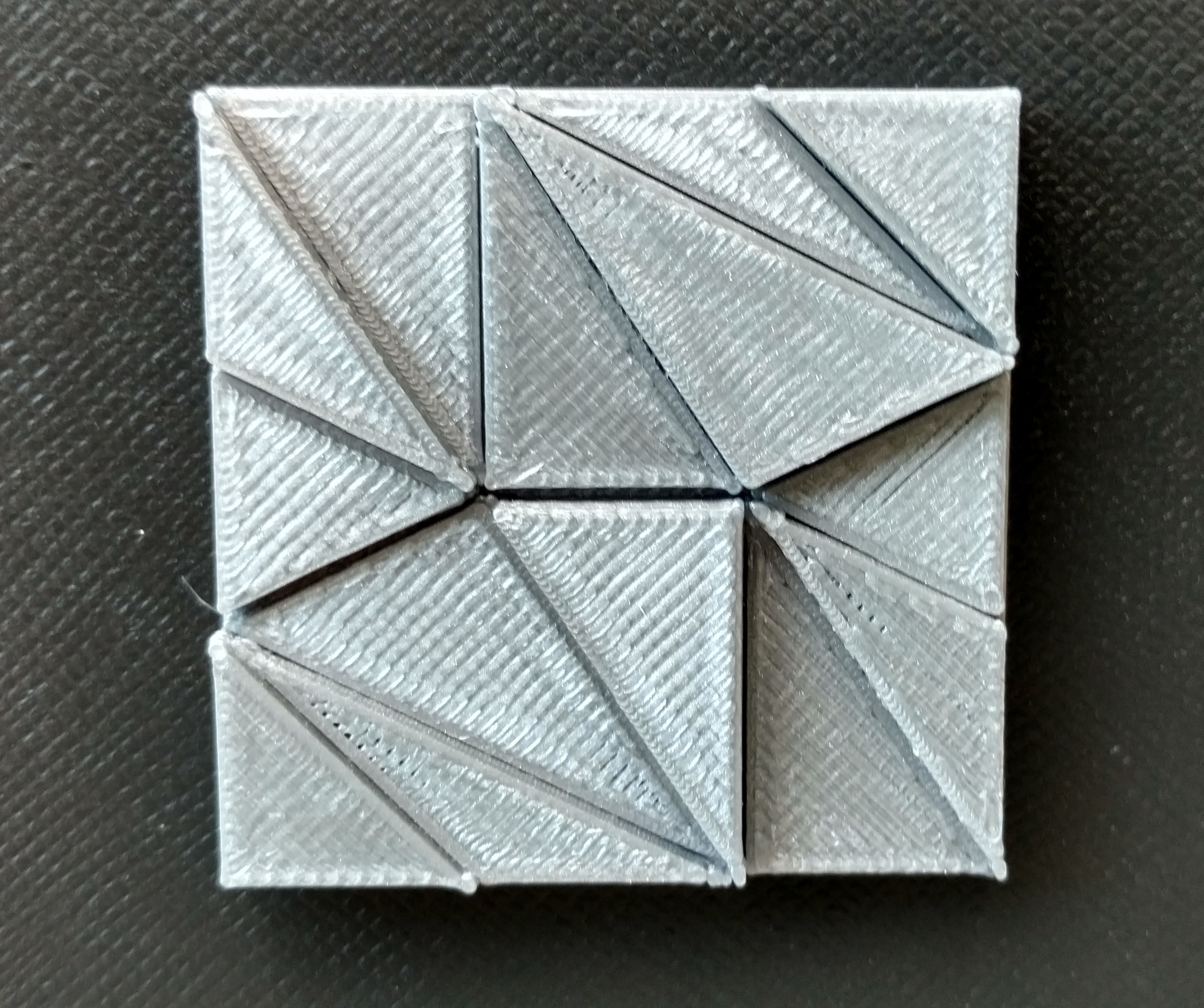}
      \caption{Grayscale becomes height.}
      \label{torus}
   \end{figure}

A JPG image can be imported into \textit{Mathematica} and converted to 0-1 grayscale, represented in a large matrix, and then this matrix, or a scalar multiple, can be used as a height function defined discretely in a table.  The \texttt{ListPlot3D} command in \textit{Mathematica} can then be used to nicely render the function.  We can do more by transforming the matrix values with some sort of filter, such as re-assigning all values greater than 0.5 with a height of 10, and all values less than 0.5 with a height of 5, using if-then-else commands.  An example of this results of this procedure can be seen in Figure \ref{torus}.

Without thinking things through, I thought I would use this approach on my son's high school portrait.  I will spare you the picture of the rendering of this function -- let's just say that the nose is not necessarily the lightest part of a portrait, the ears are not necessarily the darkest, and hair color plays extra havoc on an attempted design.  This led to the question:  how does the conversion from color to grayscale work?  And that led to the idea of \textit{luminance}.

\begin{figure}[thpb]
      \centering
      \includegraphics[scale=0.5]{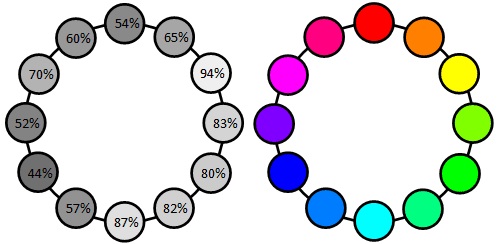}
      \caption{Colors and associated luminance. (Source: \url{http://www.workwithcolor.com/color-luminance-2233.htm}.}
      \label{color}
   \end{figure}

Photographers, painters, and other artists understand the intensity of a color is measured by its luminance. The luminance of various colors can be seen in the color wheel in Figure \ref{color}, and when these colors are projected into grayscale space, luminance becomes grayscale.  There is an excellent video by neuroscientist Margaret Livingstone on this topic \cite{Livingstone}.  Now, it would be hard to find a portrait with a lime green nose and dark blue ears.  But in the video, Livingstone mentioned a group of early 20th century artists that created the type of portraits that I could work with.

\section{The Wild Beasts of Fauvism}

\textit{Les Fauves}, or the ``Wild Beasts'' was the name given to a set of young painters, primarily in Paris, in the early twentieth century \cite{Fauvism}. They were known for their wild and unexpected use of color.  As described by Ferrier, ``The Fauves explored the spectrum; for them, colors were not only mere stimuli on the retina but could also express feelings.'' \cite{Ferrier}  The most famous Fauvist was Henri Matisse, and other well-known members of the group were Henri Manquin, Albert Marquet, Georges Braque, and Andr\'{e} Derain.  

Artists have long recognized that color and luminance play different roles in visual perception. One of the goals of Fauvism (1905-1907) was, according to Douma, ``to give color greater emotional and expressive power.'' \cite{Douma}  Douma also observes that ``Neuroscientifically, Matisse’s paintings work like a black and white photograph.''  And that is exactly what I needed to 3D print a portrait.  

   \begin{figure}[thpb]
      \centering
       \includegraphics[scale=0.25]{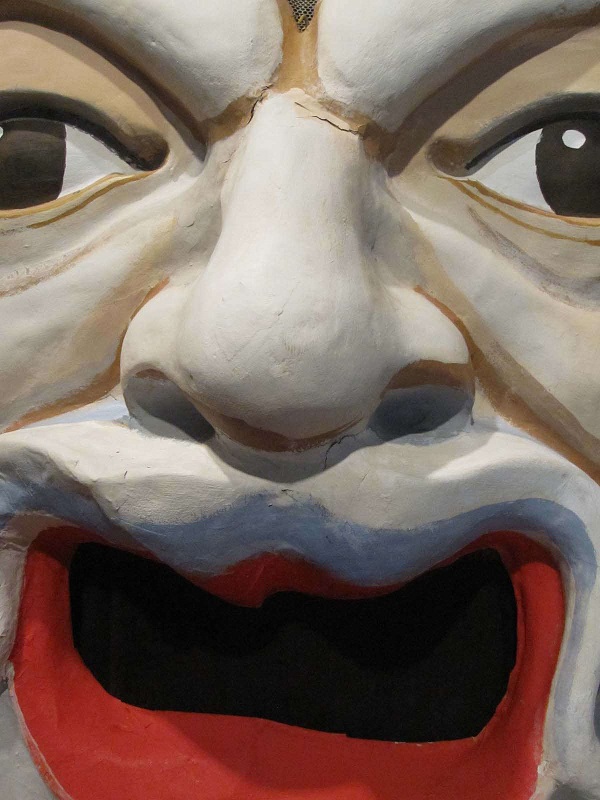}
       
       \vspace{1cm}
       
       \includegraphics[scale=0.25]{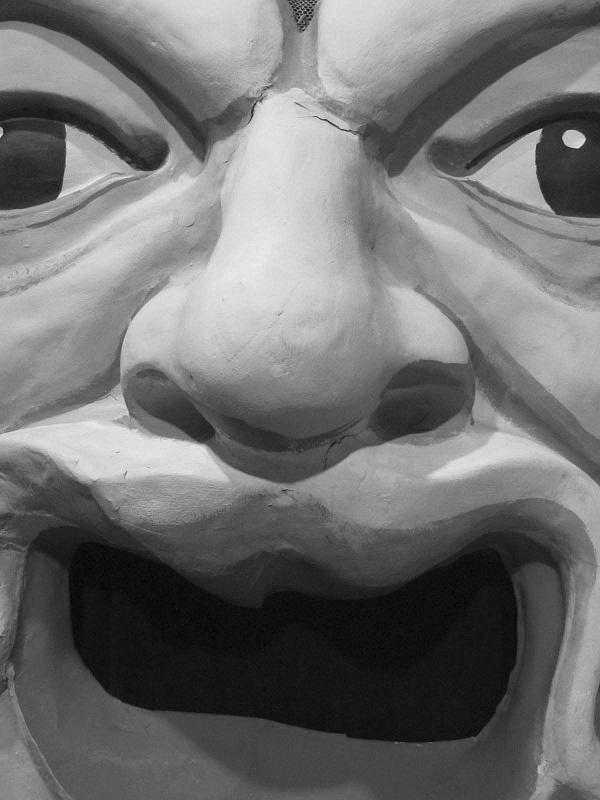}
      \caption{Original and Grayscale version of Demain's \textit{Genie Mask}.}
      \label{genie}
   \end{figure}

The portrait that I focused on for this project is \textit{Genie Mask}, seen in Figure \ref{genie}, by Andr\'{e} Derain.  Clearly, as with most Fauvist paintings, we do not have realistic colors being used for this human face.  But, the luminance of the various colors in this painting line up with pretty well with the distance from the virtual ``lens'' that is viewing the portrait.  This is apparent from the grayscale version of \textit{Genie Mask}, also in Figure \ref{genie}.

So, I was back in business, and I expected that it would not take long to find a JPG of this portrait that I could use, modify my \textit{Mathematica} file to import the image, create $f(x,y)$ and render it as a 3D image, export the result as an STL file, and 3D print the result.  For the first two steps, I was correct, but the rest of the steps required extensive computations, and actually my computer ran out of memory.  I needed a mathematical tool well-known to me to finish the project.

\section{Saving Time with Wavelets}

I've been working with wavelets since the late 1990's.  Inspired by the adoption of wavelets by the U.S. Federal Bureau of Investigation (the FBI) for compression of fingerprint images, I learned the basics using publications written by Strang \cite{Strang}, Mulcahy \cite{Mulcahy}, and eventually Daubechies \cite{Daubechies}.  After teaching introductory wavelets materials to undergraduates, Steve Schlicker and I co-wrote our own introduction to wavelets for undergraduates  \cite{DW, Ency}, and I began leading undergraduate research projects that applied wavelets in a variety of ways (see \url{https://sites.google.com/site/aboufadelreu/}).  So what are these wavelets?

\begin{figure}[thpb]
      \centering
      \includegraphics[scale=0.5]{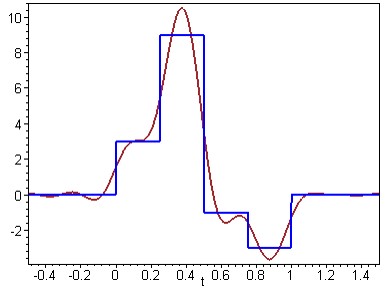}
      \caption{Example of a projection of a continuous function onto a Haar wavelet space.}
      \label{Haar}
   \end{figure}

At their simplest, a wavelet family is a set of linear independent functions that can be used for analysis, in the sense of taking other functions and writing them as a linear combination of wavelets functions.  Or, more accurately, \emph{approximating} functions with linear combinations, which is basically a projection onto a wavelet subspace.  For instance, the Haar wavelets are piecewise constant functions, and projecting continuous functions onto the Haar subspace leads to examples like what is seen in Figure \ref{Haar}.

With the type of discrete wavelet analysis being described here, the projection is created by computing a linear combination of function values.  For instance, with the Daubechies-4 wavelets, we have the following ``low-pass'' filter coefficients: 

\begin{align*}
h_0 = \frac{1+\sqrt{3}}{4 \sqrt{2}} &\approx 0.48296 \\
h_1 = \frac{3+\sqrt{3}}{4 \sqrt{2}} &\approx 0.83652   \\
h_2 = \frac{3-\sqrt{3}}{4 \sqrt{2}} &\approx 0.22414   \\
h_3 = \frac{1-\sqrt{3}}{4 \sqrt{2}} &\approx -0.1294 
\end{align*}

\noindent and low-pass filtering is calculated in this way:
\begin{equation*}
(H\mathbf{s})_k = h_0 \mathbf{s}_{2k} + h_1 \mathbf{s}_{2k+1} + h_2 \mathbf{s}_{2k+2} + h_3 \mathbf{s}_{2k+3}
\end{equation*}

There are also ``high pass'' filters which are computed by ``differencing'', as indicated in this formula from the Daubechies-4 wavelets:

\begin{equation*}
(G\mathbf{s})_k = g_{-2} \mathbf{s}_{2k-2} + g_{-1} \mathbf{s}_{2k-1} + g_0 \mathbf{s}_{2k} + g_1 \mathbf{s}_{2k+1}.
\end{equation*}
\begin{align*}
g_{-2} = \frac{1-\sqrt{3}}{4 \sqrt{2}} &\approx -0.12941  \\
g_{-1} = -\frac{3-\sqrt{3}}{4 \sqrt{2}} &\approx -0.22414  \\
g_0 = \frac{3+\sqrt{3}}{4 \sqrt{2}} &\approx  0.83652 \\
g_1 = -\frac{1+\sqrt{3}}{4 \sqrt{2}} &\approx -0.48296. 
\end{align*}

The high-pass filters are good for, among other applications, finding edges in images, and spikes in one-dimensional signals.  For instance, a high-pass filter was used in the ``Boston Pothole Project'' that I completed with three students in 2011 \cite{Boston}.

\section{Creating a 3D Version of \textit{Genie Mask}}

\begin{figure}[thpb]
      \centering
       \includegraphics[scale=0.7]{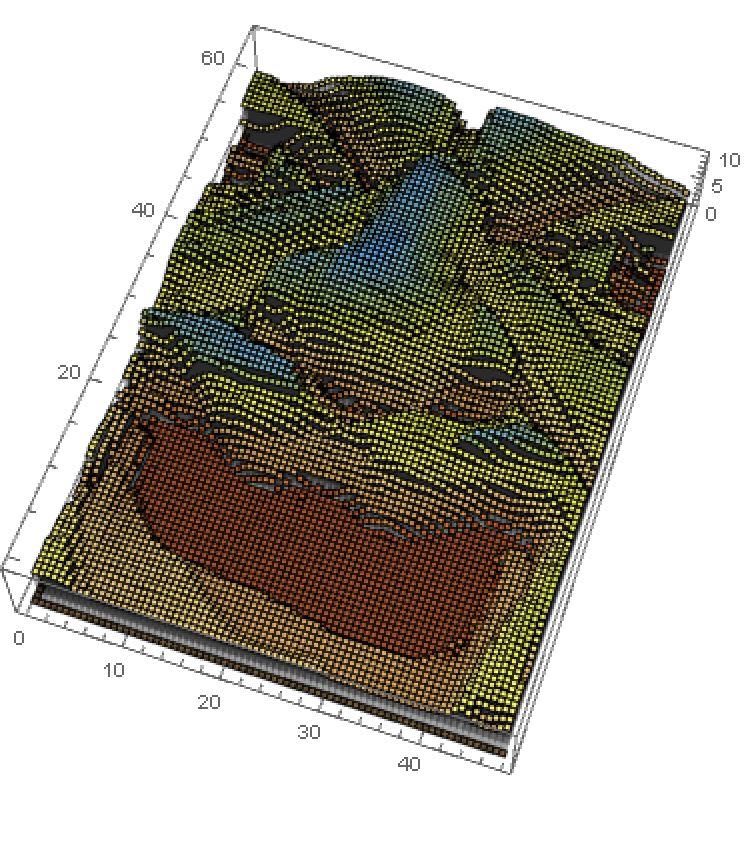}
       \vspace{0.5cm}
      \includegraphics[scale=0.7]{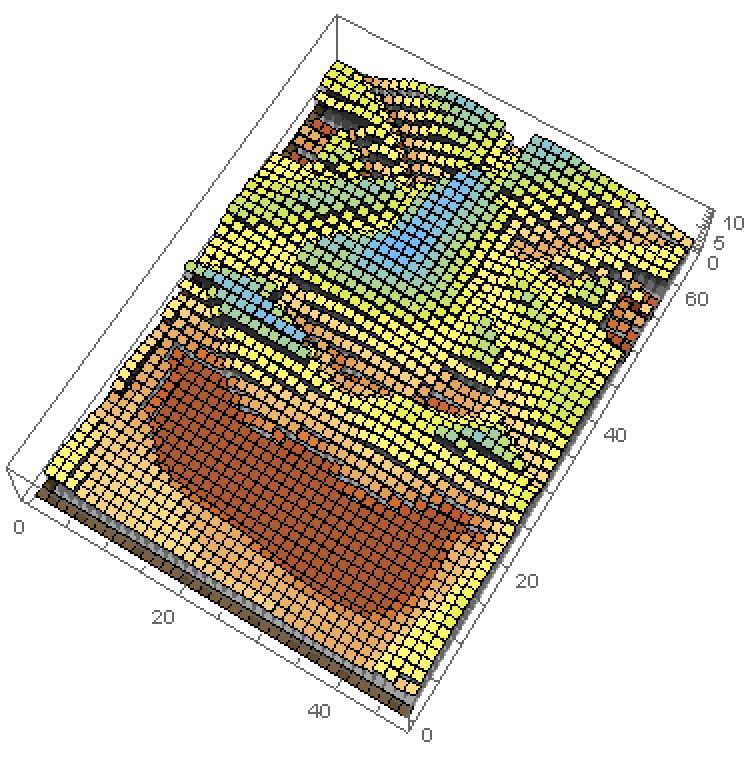}
      \caption{Rendering of output from multiple rounds of filtering of \textit{Genie Mask} matrix.}
      \label{genie2}
   \end{figure}

The filters above are for one-dimensional signals, but there are versions for higher dimensions \cite{Add}, and for my attempt to 3D print Fauvist paintings, I needed something that would work on a two-dimensional image.  For simplicity, I decided to use the two-dimensional Haar low-pass wavelets filters, which basically smooth an image by averaging pixel values on 2 by 2 blocks.  So, we think of the image as a matrix of grayscale luminances (either in the 0-1 range, or the 0-255 range, it doesn't really matter).  After dividing the matrix horizontally and vertically into 2 by 2 blocks, we replace every block with one number which is the average of the 4 values in the block.  This creates a new matrix which has half the width and height, and creates a ``smoother'' version of the matrix.  Applying this methodology several times (which is a part of what is called the ``pyramid scheme'' in the literature) leads to smoother and smoother versions of the 3D rendering of Demain's Genie Mask (see Figure \ref{genie2}).

\begin{figure}[thpb]
     \centering
      \includegraphics[scale=0.12]{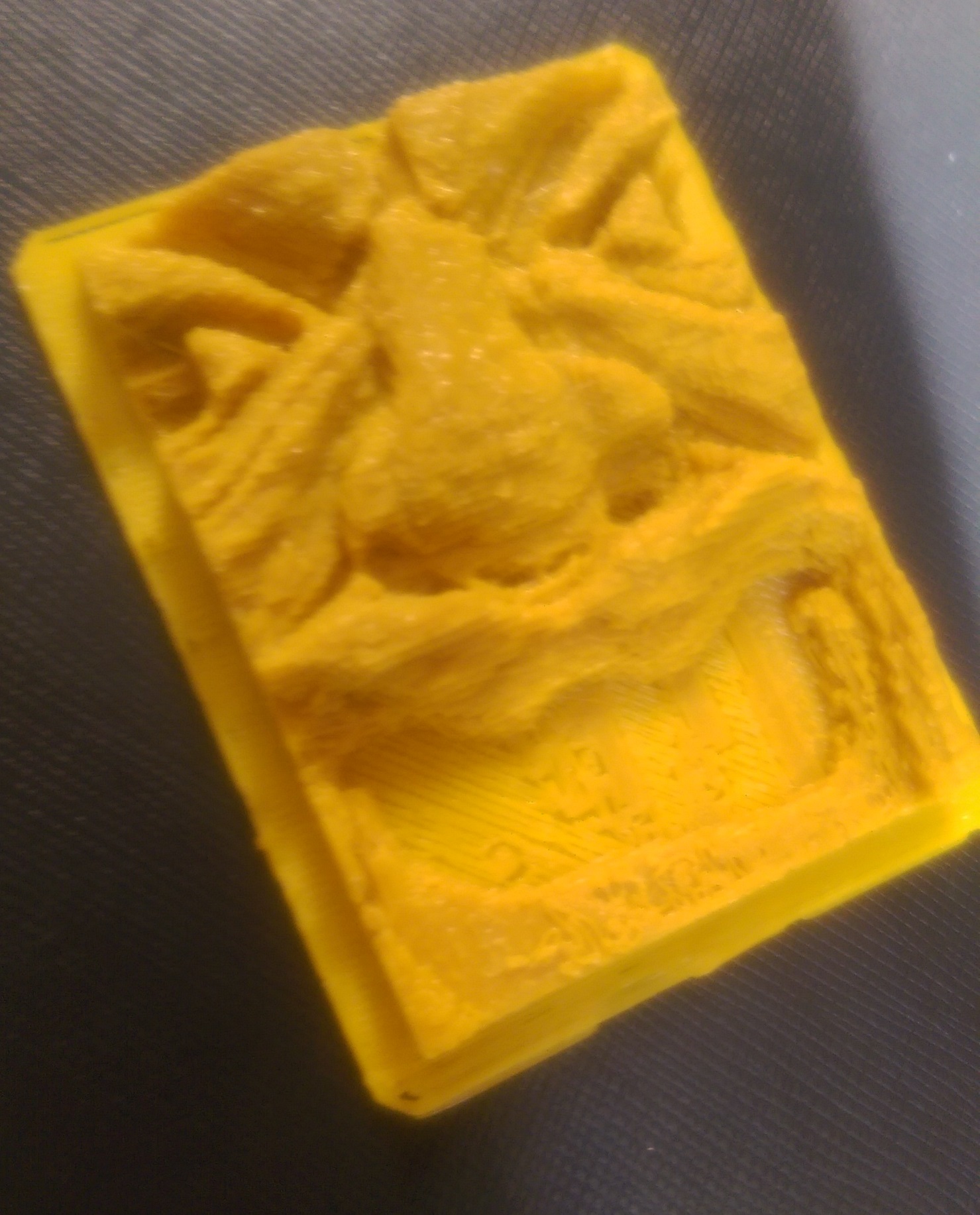}
        \caption{\textit{Genie Mask} as a 3D print.}
      \label{genie3D}
   \end{figure}

This, in turn, led to the 3D print of \emph{Genie Mask}, seen in Figure \ref{genie3D}.

\begin{figure}[thpb]
      \centering
       \includegraphics[scale=0.205]{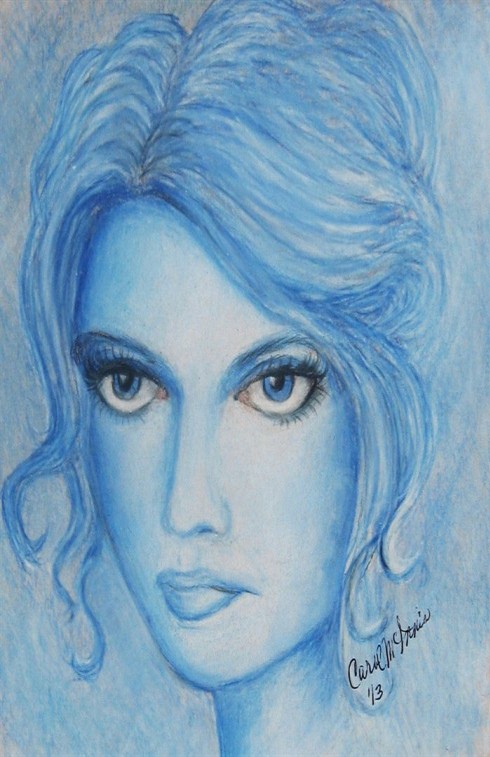}\includegraphics[scale=0.06]{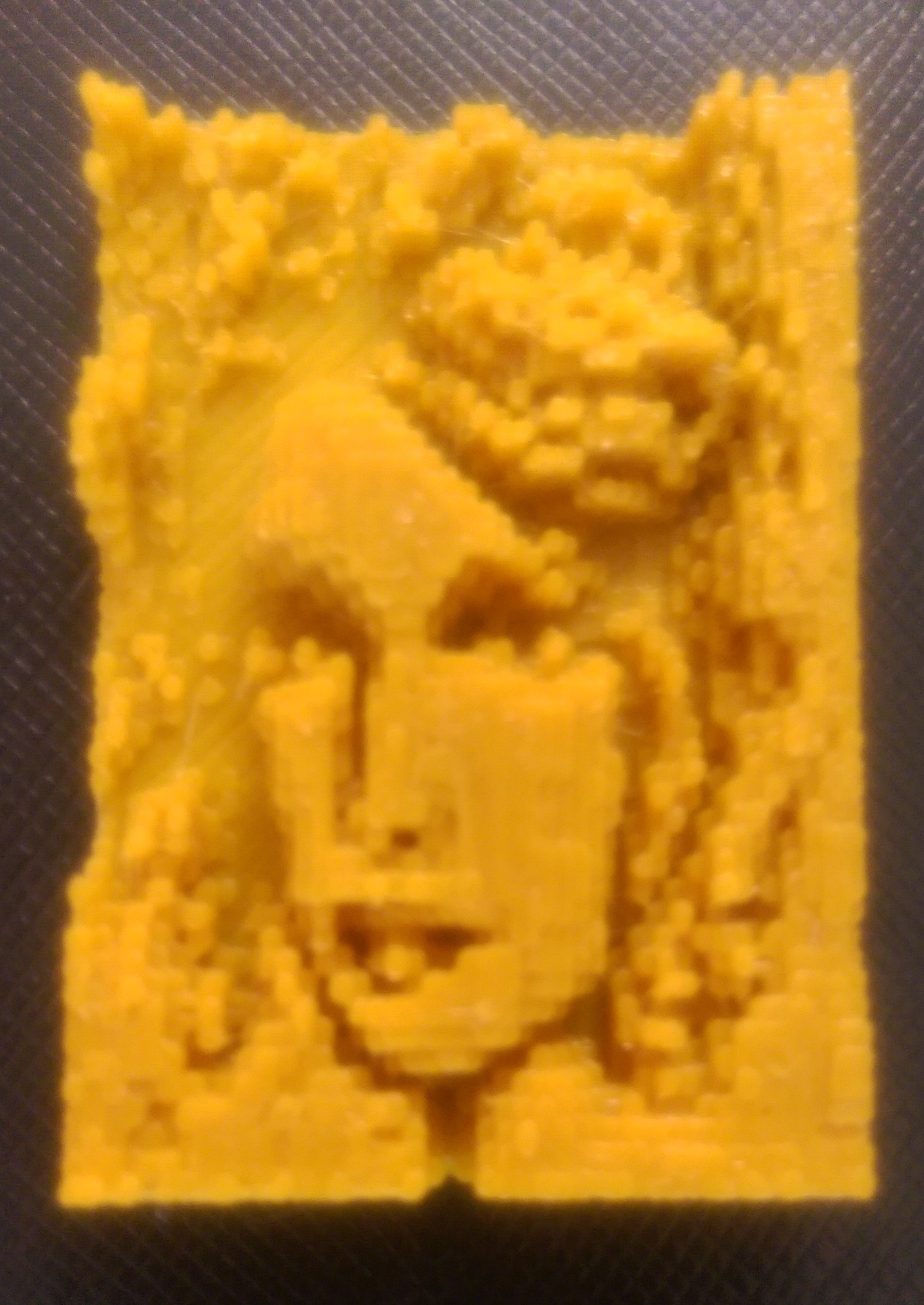}
      \caption{Another example of starting with a Fauvist image, converting to grayscale, applying wavelet filters, and then 3D printing.}
      \label{bluegirl}
   \end{figure}

In 2014, Carol McInnis, an artist that called herself ``Fauvist'', posted the ``Blue'' portrait on the left in Figure \ref{bluegirl}.  I applied my methodology to this image, yielding the 3D printed portrait also seen in Figure \ref{bluegirl}.  (Neither Carol nor this image can be found on the Internet in June 2016.  Also note that in the image that I used, I adjust the scaling of the $x$ axis.) 

\section{Final Comments}

This is not a very sophisticated application of wavelets, but it accomplished what I wanted, which is usually all that you need in a mathematical modeling project.  An interesting direction to pursue for a future project is to apply the Daubechies-4 wavelets, either the one dimensional in two directions or the two-dimensional version \cite{biva}, to see how that would smooth the original image. One advantage of using these wavelets is that they are applied to 4 by 4 boxes rather than 2 by 2, so that fewer rounds of smoothing would be needed.

Some other directions to pursue:  Is there a good ``Fauvist filter'' to apply to school portraits?  If not, is it possible to create one?  And, of course, would this approach work with other Fauvist paintings, particularly those which are not portraits?


\end{document}